\documentclass[12pt]{article}
\usepackage{latexsym}
   \let\ud=\d       
\let\a=\alpha
\let\b=\beta

\let\d=\delta       \let\vd=\partial             \let\D=\Delta

\let\lm=\lambda                                  
\let\m=\mu

\let\p=\pi

\let\s=\sigma                   

\let\w=\omega

\def\cB{{\cal B}}
\def\rd{{\rm d}}
\def\cF{{\cal F}}
\def\ie{{\it i.e.}}
\def\sh{s^\sharp}
\def\cM{{\cal M}}
\def\shM{{\cal M}^\sharp}
\def\nM{{\cal M}^\natural}
\def\cV{{\cal V}}

\def\shx{x^\sharp}
\def\sll#1{\rlap{\,\raise1pt\hbox{/}}{#1}}
\def\inv#1{{\ttt{1\over#1}}}

\def\beq#1{\begin{equation}\label{#1}}
\def\eeq{\end{equation}}

\let\To=\Rightarrow

\let\into=\hookrightarrow

\def\dual{\hbox{\Large$*$}}
\def\define{\buildrel\rm def\over =}
\def\dual{\buildrel*\over\sim}
\def\Dual{\buildrel\star\over\sim}
\def\vev#1{\left\langle#1\right\rangle}
\let\sss=\scriptscriptstyle
\let\SSS=\scriptstyle
\let\ttt=\textstyle

\let\ba=\overline
\let\2=\underline
\let\id=\equiv

\let\Tw=\widetilde
\let\3=\ldots
\def\rlx{\relax\leavevmode}
\def\inbar{\vrule height1.5ex width.8pt depth0pt}
\def\sinbar{\vrule height1ex width.45pt depth0pt}
\def\ssinbar{\vrule height.7ex width.35pt depth0pt}
\def\IC{{\relax\leavevmode
                \ifmmode\mathchoice
                       {\hbox{\kern.25em\inbar\kern-.3em{\sf C}}}
                       {\hbox{\kern.25em\inbar\kern-.3em{\sf C}}}
                       {\hbox{\kern.20em\sinbar\kern-.25em$\SSS\sf C$}}
                       {\hbox{\kern.18em\ssinbar\kern-.22em$\sss\sf 
C$}}
                \else{\hbox{\kern.3em\inbar\kern-.3em{\rm C}}}\fi}}
\def\Ik{{\rlx{\rm I\kern-.18em k}}}
\def\lbar{{\mathaccent"7020\lambda}}
\def\IP{{\rlx{\rm I\kern-.18em P}}}
\def\IQ{{\relax\leavevmode
                \ifmmode\mathchoice
                       {\hbox{\kern.33em\inbar\kern-.3em{\rm Q}}}
                       {\hbox{\kern.33em\inbar\kern-.3em{\rm Q}}}
                       {\hbox{\kern.28em\sinbar\kern-.25em$\SSS\rm Q$}}
                       {\hbox{\kern.25em\ssinbar\kern-.22em$\sss\rm 
Q$}}
                \else{\hbox{\kern.3em\inbar\kern-.3em{\rm C}}}\fi}}
\def\IR{{\rlx{\rm I\kern-.18em R}}}
\def\ZZ{{\relax\leavevmode
                \ifmmode\mathchoice
                       {\hbox{\sf Z\kern-.4em Z}}
                       {\hbox{\sf Z\kern-.4em Z}}
                       {\lower.9pt\hbox{\scriptsize\sf Z\kern-.36em Z}}
                       {\lower1.2pt\hbox{\tiny\sf Z\kern-.36em Z}}
                \else{\sf Z\kern-.4em Z}\fi}}
\def\Ione{{\rlx{\rm 1\kern-3pt l}}}

\def\NP#1{{\itshape Nucl.\,Phys.\,}{\bf#1\,}}

\def\MPL#1{{\itshape Mod.\,Phys.\,Lett.\,}{\bf#1\,}}

\def\CP#1#2{\relax\ifmmode\IP^{#1}_{#2}\else\IP$^{#1}_{#2}$\fi}
\def\topic#1{\par\noindent\underline{#1}\nobreak\vglue0mm%
               \noindent\ignorespaces}
\def\ytem#1{\par\noindent\hbox to\parindent{\hss\em#1.~}\ignorespaces}
\def\yytem#1{\par\noindent\hbox to2\parindent{\hss\em#1.~}\ignorespaces}
\newtheorem{defn}{Definition}

\newtheorem{lemma}{Lemma}
\newtheorem{thrm}{Theorem}
\newtheorem{corl}{Corollary}
\font\msam = msam10 scaled 1200
  \let\chck\undefined  \def\chck{{\msam\char"58}}
\def\QED{\relax\leavevmode\mbox{}\hfill
          $\raisebox{2.5pt}{\chck}\mkern-18mu\Box$}
\def\Proof{\par\noindent{\bf Proof:~}\ignorespaces}
\def\Remark{\par\noindent{\bf Remark:~}\ignorespaces}
\thicklines     
\def\Label#1{\label{#1}%
    \smash{\hbox to0pt{\raise1ex\hbox{\tiny\{#1\}}\hss}}}
\def\noLabels{\let\Label=\label}

\def\Eq#1{Eq.~(\ref{#1})}
\def\Eqs#1{Eqs.~(\ref{#1})}

\catcode`@=11   
\newbox\t@b@x
\def\rightarrowfill{$\m@th \mathord- \mkern-6mu
      \cleaders\hbox{$\mkern-2mu \mathord- \mkern-2mu$}\hfill
       \mkern-6mu \mathord\rightarrow$}
\def\tooo#1{\setbox\t@b@x=\hbox{\footnotesize$#1$}%
              \mathrel{\mathop{\hbox to\wd\t@b@x{\rightarrowfill}}%
               \limits^{#1}}\,}
\def\leftarrowfill{$\m@th \mathord\leftarrow \mkern-6mu
      \cleaders\hbox{$\mkern-2mu \mathord- \mkern-2mu$}\hfill
       \mkern-6mu \mathord-$}
\def\froo#1{\setbox\t@b@x=\hbox{\footnotesize$#1$}%
              \mathrel{\mathop{\hbox to\wd\t@b@x{\leftarrowfill}}%
               \limits^{#1}}\,}
\catcode`@=12                   

\noLabels
\input epsf.tex
\begin{document}
\leftline{math.AG/0210394}

\begin{center}

{\LARGE\bf On the Geometry and Homology of\\[2mm]
            Certain Simple Stratified Varieties}\\[10mm]

{\bf T. H\"{u}bsch}\footnote{thubsch@howard.edu; On leave from the
    ``Rudjer Bo\v skovi\'c'' Institute, Zagreb, Croatia.}
      and
{\bf A.~Ra\ud{h}m\={a}n}\footnote{arahman@howard.edu}\\[1mm]
{\it  Department of Mathematics and\\
      Department of Physics and Astronomy\\
      Howard University\\
      Washington, DC 20059}\\[5mm]

{\bf ABSTRACT}\\[3mm]
\parbox{5in}{
 We study certain mild degenerations of algebraic varieties which appear
in  the analysis of a large class of supersymmetric theories, including
superstring theory. We analyze Witten's $\s$-model~\cite{wit:phases} and
find that the non-transversality of the superpotential induces a
singularization and stratification of the ground state variety.
 This stratified variety (the union of the singular ground state variety
and its exo-curve strata) admit homology groups which, excepting the
middle dimension, satisfy the ``K\"ahler package'' of
requirements~\cite{gmp:strata}, extend the ``flopped'' pair of small
resolutions to an ``(exo)flopped'' triple, and is compatible with mirror
symmetry~\cite{yau:mirror} and string theory~\cite{str:nullBH,hub:strcoh}.
 Finally, we revisit the conifold transition \cite{cgh:confld} as
it applies to our formalism.

}
\end{center}

\bibliographystyle{plain}

\section{Introduction, Results and Summary}
\Label{s:IRS}
In string theory, rather than being an assumed arena, the spacetime is
identified with the dynamically determined `ground state variety' of
a (supersymmetric) $\s$-model~\cite{gsw:string,wit:phases,pol:string}.
  In the simplest physically interesting and nontrivial
case~\cite{chsw:strcys,wit:phases}, the spacetime is of the form
$M^{3,1}{\times}K$, where $K$ is a compact Calabi-Yau 3-fold modeled
from the (bosonic subset of the) `field space' of the
$\s$-model\footnote{To avoid obscuringly complicated notation,
we focus on a simple example and discuss generalizations later.},
$\cF=\{p,s_0,{\cdots},s_4\}\simeq\IC^6$, which admits a $\IC^*$
action:
\begin{equation}
   \hat{\lm}:~\{p,s_0,{\cdots},s_4\} \mapsto
           \{\lm^{-5}p,\lm s_0,{\cdots},\lm s_5\}~,\qquad \lm\in\IC^*~.
   \Label{e:ProjMap}
\end{equation}
The `ground state variety' is defined to be~\cite{wit:phases,agm:chngtp}
\begin{equation}
   \cV~\define~[\,(\vd W)^{-1}(0)-0\,]/\hat{\lm}~,
   \Label{e:HolQ}
\end{equation}
with the $\hat{\lm}$-invariant holomorphic `superpotential'
\begin{equation}
   W~\define~p{\cdot}G(s)~. \Label{e:W=pG}
\end{equation}
Alternatively, we denote by $\hat{\lbar}$ the $|\lm|=1$ restriction of
the map~(\ref{e:ProjMap}), and define the `potential'
\begin{equation}
 U_r~\define~\|\vd W\|^2+D_r^2~, \Label{e:U=dW+D}
\end{equation}
where
\begin{equation}
   D_r~\define~\|s\|^2-5|p|^2-r~,\qquad r\in\IR~. \Label{e:gaugeD}
\end{equation}
Then
\begin{equation}
   \cV~\simeq~[\,U_r^{-1}(0)-0\,]/\hat{\lbar}~.
   \Label{e:SimQ}
\end{equation}

Owing to the positive definiteness of $U_r$,
\begin{equation}
   U_r^{-1}(0) ~=~ (\vd W)^{-1}(0) ~\cap~ D_r^{-1}(0)~. \Label{e:RMfld}
\end{equation}
Furthermore, the $\hat{\lm}$-invariance of $W\,{=}\,pG$ implies that
$G(s)$ is a degree-5 homogeneous complex polynomial
\begin{equation}
   G(\lm s_0,{\cdots},\lm s_4) ~=~ \lm^5 G(s_0,{\cdots},s_4)~,
   \Label{e:HomHoloG}
\end{equation}
whereupon the zero locus of $\vd W$ is the intersection of the cones
\begin{equation}
   (\vd W)^{-1}(0)~=~G^{-1}(0)~\cap~(p{\cdot}\vd_sG)^{-1}(0)~.
\end{equation}
The above definition may then be rephrased as follows:
\begin{defn}\Label{d:GndStVar}
   Given the polynomials $G(s)$ and $D_r$ as defined in
Eqs.~(\ref{e:HomHoloG}) and~(\ref{e:gaugeD}), respectively,
   the `ground state variety' is
\begin{equation}
  \begin{array}{rcl}
   \cV&=&\Big\{G^{-1}(0)\cap
               (p{\cdot}\vd_sG)^{-1}(0)-0\,\Big\}\Big/\hat{\lm}~,\\[2mm]
      &=&\Big\{G^{-1}(0)\cap
               (p{\cdot}\vd_sG)^{-1}(0)\cap
               D_r^{-1}(0)-0\,\Big\}\Big/\hat{\lbar}~,
  \Label{e:GndStVar}
  \end{array}
\end{equation}
where the $S^1$-action, $\hat{\lbar}$, in the latter (symplectic) quotient
is the $|\lm|=1$ restriction of the
$\IC^*$-action~(\ref{e:ProjMap}) in the former (holomorphic) quotient.

$\cV^+$ ($\cV^-$) shall denote the restriction of $\cV$ to positive (negative) values of $r$ in \Eq{e:gaugeD}.
\end{defn}

\Remark
As we show in more detail in Section~\ref{s:GSV}, the dependence $U_r$,
in \Eq{e:gaugeD}, turns the `ground state variety' into a 1-parameter family
of (stratified) varieties\footnote{The general category of `stratified
varieties' is specified for example in the works~\cite{gmp:strata}.
Our situation is far simpler: we will encounter unions of several
(complex, algebraic) varieties of complex dimension $0,\cdots,3$,
possibly connected at codimension\,$\geq1$ subspaces.}, and the
subtraction of zero in Eqs.~(\ref{e:GndStVar}) separates the two
branches, $\cV^\pm$, defined with $r{>}0$ and $r{<}0$, respectively.
Moreover, the dependence on $r$, as defined originally in the gauged
linear $\s$-model~\cite{wit:phases}, is complicated near $r{=}0$
by quantum corrections and we restrict to $r{\neq}0$.

Generalizations involve ({\it a})~additional $p$ and $s$ variables,
({\it b})~additional corresponding terms in the superpotential~(\ref{e:W=pG}),
and ({\it c})~additional maps~(\ref{e:ProjMap}) and their modifications where
the exponents of $\lm$ are different integers: negative for the $p$'s,
positive for the $s$'s.
The generalization~({\it a}) turns the $\IC^*$ action~(\ref{e:ProjMap})
into a more general toric action, while ({\it a})~enlarges the field space
and ({\it b})~modifies the ``moment map''~(\ref{e:W=pG}) accordingly.
The resulting `ground state varieties' will thus include intersections of
hypersurfaces in products of toric
varieties~\cite{wit:phases,agm:chngtp,ful:torvar,gun:torvar}.
All of  these are of the form given in the definition~\ref{d:GndStVar},
with the `ingredients' $\{(p,s),\hat\lm,D_r,W\}$ duly modified.

Our main result is:
\begin{thrm}\Label{t:main}
Let $\cV$, the `ground state variety' of the gauged linear
$\s$-model~\cite{wit:phases}, be defined in Definition~\ref{d:GndStVar}.
Then
  \ytem{1} $\cV^+$ is a stratified variety~\cite{gmp:strata} when the
polynomial $G(s)$ is non-transversal at $n$ isolated rays, $\sh$.
  \ytem{2} Then, $\cV^+=\shM\cup\bigcup_iA_i$ with $\shM=G^{-1}(0)$
smooth except at $n$ isolated nodes, where the $n$ noncompact `antennae'
$A_i\simeq\IC^1$ attach.
  \ytem{3} For $\dim_\IC\cV^+{=}3$, the minimal holomorphic
compactification $\ba{\cV}^+=\shM\cup\bigcup_i\bar{A}_i$ satisfies:
  \yytem{a} $\ba{\cV}^+$ is an exoflop of the small resolution(s) of $\shM$
in the sense of Ref.~\cite{agm:chngtp}, and
  \yytem{b} $\oplus_qH^{2q}(\ba{\cV}^+)$ satisfies the
``K\"ahler package'' of requirements~\cite{gmp:strata},
and is compatible with mirror symmetry~\cite{yau:mirror}
and string theory~\cite{str:nullBH,hub:strcoh}.
\end{thrm}

This article is organized as follows: Section~\ref{s:GSV} shows that the
`ground state variety' becomes stratified as $G(s)$ becomes
non-transversal, and we explore the induced (exo-)strata\footnote{We will
use the prefix `exo' to denote (components of) strata that are `external'
to the `main' stratum.} and their union. Section~\ref{s:Hom} explores the
contribution of the induced (exo-)strata to the homology of the `ground
state variety.' Section~\ref{s:DefRes} re-examines the `conifold
transition' of Refs.~\cite{cgh:confld} in view of Theorem~\ref{t:main}.

\section{The Ground State Variety}
\Label{s:GSV}
We now turn to analyze the geometry of the ground state variety, as
determined by the choice of the homogeneous holomorphic polynomial 
$G(s)$.
Such polynomials typically depend on a multitude of parameters; when 
properly
accounted for redundancies, these span (a subspace of) the moduli space 
of
the ground state variety. Thus, we automatically have a family of 
ground
state varieties, fibered over this (partial) moduli spaces. Works in 
the
literature, Ref.~\cite{wit:phases} and the subsequent studies, all 
assumed
$G(s)$ to be transversal and so have explored the generic fibre of this
family. We begin by analyzing this case in some detail, and then turn 
to
the less generic mild degenerations of the fibre.

\subsection{The transversal case}
\Label{s:TransV}
$G(s)$ being transversal\footnote{Transversality ensures
that the projective hypersurface defined by $G=0$ is smooth.}, 
$G= \rd G=0$ only at $s=0$.\ In this case,
the zero locus of $\vd W=(G,p{\cdot}\vd_s G)$ is a union of two
branches\footnote{The subsequent analysis for non-transversal $G$,
the case of our real interest, is more detailed and shown below.
The Reader can then recover the presently omitted details as 
a special case; see also Ref.~\cite{wit:phases}.}:
\begin{equation}
(G)^{-1}(0)\cap(p{\cdot}\vd_sG)^{-1}(0)
  ~=~\Big\{p{=}0,~s:G(s)=0\Big\}\cup\Big\{p,~s{=}0\Big\}~.
  \Label{e:union}
\end{equation}

So, following the first (holomorphic quotient) part of
definition~\ref{d:GndStVar}, we have that
\begin{equation}
   \cV~=~\Big\{p{=}0,~s:G(s)=0\Big\}\Big/\hat{\lm}
          ~\cup~\Big\{p,~s{=}0\Big\}\Big/\hat{\lm}~,
   \Label{e:qunion}
\end{equation}
where the quotients are taken after the fixed point of the
$\hat{\lm}$-action, $\{s,p=0\}$, is excised.
  Now, since $\{s{\neq}0\}/\hat{\lm}$ is $\IP^4$, then
\begin{equation}
   \Big\{p{=}0,~s\neq0:G(s)=0\Big\}\Big/\hat{\lm}~=~\cM \Label{e:CYM}
\end{equation}
is the Calabi-Yau quintic hypersurface in $\IP^4$.

The second term in the union~(\ref{e:qunion}) is
\begin{equation}
   \Big\{p{\neq}0,~s{=}0\Big\}\Big/\hat{\lm}~\simeq~\IC^*/\IC^*
    ~\simeq~\{\hbox{pt.}\}/\ZZ_5~. \Label{e:LGO}
\end{equation}
This is the `fuzzy point'~\cite{agm:chngtp} of the Landau-Ginzburg 
orbifold.
$\ZZ_5$ is the subroup of $\hat\lm$ which leaves $G(s)$
invariant and so acts trivially on both $W{=}p{\cdot}G(s)$ and on $p$.
\begin{figure}[ht]
  \setlength{\unitlength}{.85mm}
  \begin{picture}(170,125)(-20,0)
        \put(0,0){\epsfxsize=140mm\epsfbox{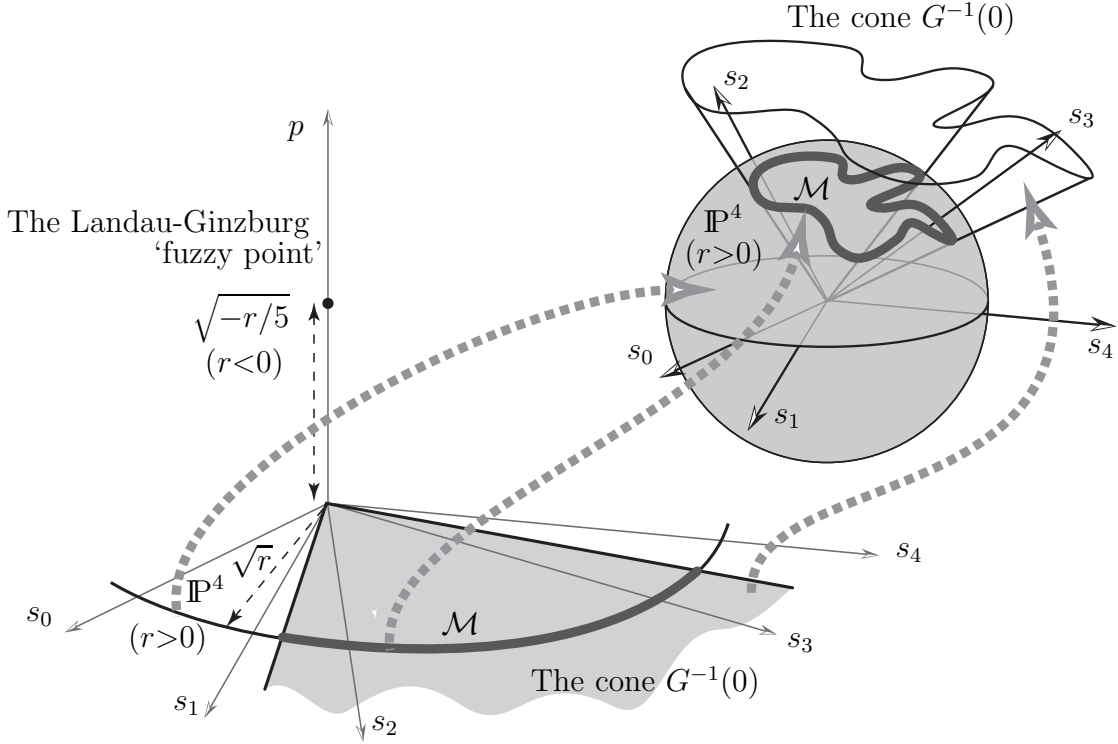}}
        \put(-6,19){$s_0$}
        \put(17,5){$s_1$}
        \put(48,2){$s_2$}
        \put(113,15){$s_3$}
        \put(130,29){$s_4$}
        \put(73,8){The cone $G^{-1}(0)$}
        \put(19,22){$\IP^4$}
        \put(10,15){$(r{>}0)$}
        \put(26,27){$\sqrt{r}$}
        \put(59,17){$\cM$}
        \put(35,95){$p$}
        \put(-9,80){The Landau-Ginzburg}
        \put(14,75){`fuzzy point'}
        \put(20,65){$\sqrt{-r/5}$}
        \put(22,58){$(r{<}0)$}
        \put(88,60){$s_0$}
        \put(111,50){$s_1$}
        \put(103,103){$s_2$}
        \put(157,97){$s_3$}
        \put(160,61){$s_4$}
        \put(113,112){The cone $G^{-1}(0)$}
        \put(100,80){$\IP^4$}
        \put(97,75){$(r{>}0)$}
        \put(114,85){$\cM$}
  \end{picture}
  \caption{The ground state variety, $\cV\in{\cal F}$, and its
`geometric' phase, $\cM\in{\cal F}|_{p=0}$, in the top right inset.
This inset represents $\cM=G^{-1}(0)\cap\IP^4$, in
${\cal F}|_{p=0}\approx\IC^5$ spanned by $s=(s_0,\cdots,s_4)$.
The dashed grey arrows identify the image under the
projection along $p$ of this in the full field space,
$\cF\approx\IC^6$, spanned by $(p,s_0,{\ldots},s_5)$.}
\Label{f:smooth}
\end{figure}

The two quotients in the union~(\ref{e:qunion}) are thus manifestly
disconnected: the former, (\ref{e:CYM}), lies entirely in the
$\{p{=}0,~s{\neq}0\}$-subspace of the field space $\cF$, whereas the
latter, (\ref{e:LGO}), lies well in the complementary
$\{p{\neq}0,~s{=}0\}$-subspace.
The above is illustrated in Fig.~\ref{f:smooth}. Even in the 
transversal
case, $\cV$ may be regarded as a stratified variety, consisting of two
disconnected objects: a (complex) 3-dimensional one and a (complex)
0-dimensional one, each of which containing a single variety: $\cM$ and
$\{\hbox{pt.}\}/\ZZ_5$, repsectively.

In fact, the second component, $\{\hbox{pt.}\}/\ZZ_5$, actually lies in 
the
`second sheet' of the field space $\cF$. To see this, it will be useful 
to
also present $\cV$ using the alternate (symplectic quotient)
definition~\ref{d:GndStVar}:
\begin{enumerate}
  \item When $r\gg0$, $D_r^{-1}(0)\neq0$ implies that $\|s\|^2\neq0$, 
and
so $\vd_sG\neq0$ as $G$ is transversal. Then, $(\vd W)^{-1}(0)$ lies
entirely in the $(p{=}0)$ $s$-hyperplane, and $\cV$ is the
$\hat\lm$-quotient, \ie, the complex base of the cone $G^{-1}(0)$.
Then, $\cV=\cM\id[\{p=0\}\cap G^{-1}(0)]/\hat{\lm}$. Since
$\{s{\neq}0\}/\hat{\lm}=\IP^4$, the  {\it projective\/} Calabi-Yau 
quintic
hypersurface~(\ref{e:CYM}) is $G^{-1}(0)/\hat{\lm}=\cM\into\IP^4$.
  \item When $r\ll0$, $D_r^{-1}(0)\neq0$ implies that $|p|^2\neq0$, 
and so
$\|s\|^2{=}0$ since $G$ is transversal. Then, $(\vd W)^{-1}(0)$ lies
entirely in complex the $p$-plane, and is the
`fuzzy point'~\cite{agm:chngtp}, $\Big\{|p|=\sqrt{|r|/5}\Big\}/\ZZ_5$,
of the Landau-Ginzburg orbifold~(\ref{e:LGO}).
\end{enumerate}

At the critical point $r=0$, these two branches formally collapse to 
the
the highly degenerate point $p,s=0$, which is the branching point of 
the two
`sheets' of the field space $\cF$. This point is, by definition, 
excised
before taking the quotients~(\ref{e:HolQ}). Indeed, for 
applications to
string theory, the preceding analysis is not to be trusted in the 
region
near $p,s{=}0$ since quantum corrections modify the 
map~(\ref{e:gaugeD}) and
so also the structure of the quotients in definition~\ref{d:GndStVar}; 
see
Ref.~\cite{wit:phases}. For this reason, we will mostly concern
ourselves with the $r\gg0$ `sheet', and comment on occasion on the
$r\ll0$ `sheet', but leave any `connection' between the two `sheets'
unexplored for now.

\subsection{The conifold with exocurves}
\Label{s:SingV}
Unlike Ref.~\cite{wit:phases} and subsequent work, we will be concerned
with ground state varieties using homogeneous holomorphic polynomials 
$G(s)$
which are non-transversal along $n$ isolated (complex) directions:
\begin{equation}
   \vd G(s)~=~0 \qquad\To\qquad s~=~\sh_j~,\quad j=1,{\cdots},n~.
   \Label{e:nodalG}
\end{equation}
Clearly, $\{\sh_j\}\simeq\IC^1$, and we denote by
$\cB^n\define\sqcup_{j=1}^n\{\sh_j\}$ the `bouquet' of $n$ $\IC^1$'s 
all
meeting at the origin. Since $G(s)$ is holomorphic and homogeneous,
$\vd G(s)=0$ implies $s{\cdot}\vd_sG(s)=5G(s)=0$ and the $G(s)=0$
condition is automatically satisfied on $\cB^n$. Thus, we find that
\begin{equation}
   (G)^{-1}(0)\cap(p{\cdot}\vd_sG)^{-1}(0)  
~=~\Big\{p{=}0,~s\neq\sh_j:G(s)=0\Big\}\cup\Big\{\{p\}{\times}\cB^n\Big\}~.
  \Label{e:Union}
\end{equation}

So, following the first (holomorphic quotient) part of
definition~\ref{d:GndStVar}, we now have that
\begin{equation}
   \cV~=~\Big\{p{=}0,~s\neq\sh_j:G(s)=0\Big\}\Big/\hat{\lm}
          ~\cup~\Big\{\{p\}{\times}\cB^n\Big\}\Big/\hat{\lm}~;
   \Label{e:QUnion}
\end{equation}
again, the quotients are taken after the fixed point of the
$\hat{\lm}$-action, $\{s,p=0\}$, is excised.
  Now, since $\{s{\neq}0,\sh_j\}/\hat{\lm}$ equals $\IP^4$ without
its points where $G(s)|_{\IP^4}$ is non-transversal, then
\begin{equation}
   \Big\{p{=}0,~s\neq0,\sh_j:G(s)=0\Big\}\Big/\hat{\lm}
   ~=~\shM-\hbox{Sing}(\shM) \Label{e:nonConifold}
\end{equation}
is the non-singular (and non-compact) part of the
conifold\footnote{Following Ref.~\cite{cgh:confld}, a {\it conifold\/} 
is a
variety which is smooth except for a finite number of isolated conical
singularities. Furthermore, herein we will consider only varieties with
nodes (double points).} $\shM\into\IP^4$. Note that
\begin{equation}
  \begin{array}{rcl}
   \hbox{Sing}(\shM)
    &\define& \Big\{p{=}0,~s{=}\sh_j:G(s)=0\Big\}\Big/\hat{\lm}\\[2mm]
    &=& \bigcup_{j=1}^n\Big\{p{=}0,~s{=}\sh_j\Big\}\Big/\hat{\lm}
     ~=~ \bigcup_{j=1}^n \shx_j \subset \shM~, \Label{e:SingM}
  \end{array}
\end{equation}
since $G(\sh_j){=}0$; $\shx_j$ are the singular points of
\begin{equation}
   \shM~\define~[\{p=0\}\cap G^{-1}(0)]/\hat{\lm}~. \Label{e:shM}
\end{equation}

The second quotient in the union~(\ref{e:QUnion}) is quite more 
intricate. Setting $p{=}0$ in
\Eq{e:gaugeD}, we see that $\{\{p\}{\times}\cB^n\}/\hat{\lm}$ is
non-empty in the $r{>}0$ `sheet' of the field space $\cF$, and also 
that it
includes the points 
$\{p{=}0,\sh_j{\neq}0\}/\hat{\lm}=\hbox{Sing}(\shM)$. On
the other hand,
\begin{equation}
   \Big(\Big\{\{p\}{\times}\cB^n\Big\}\Big/\hat{\lm}\Big)
   ~=~\sqcup_{j=1}^n A_j~,\qquad
    A_j~\define~\{p,\sh_j\}\Big/\hat{\lm}~. \Label{e:Antennae}
\end{equation}
Each of the $A_j$'s contains precisely one of the singular points of 
$\shM$,
as given in \Eq{e:SingM}
\begin{equation}
   \shx_j~=~\{p{=}0,~\sh_j\}\Big/\hat{\lm}~=~A_j\cap\shM~.
\end{equation}
Thus, ground state
variety~(\ref{e:QUnion}), which is the {\it connected\/} union
of~(\ref{e:Antennae}) and of~(\ref{e:nonConifold}), is then
\begin{equation}
   \cV~=~\shM ~\cup~ \sqcup_{j=1}^n A_j~. \Label{e:AntennaM}
\end{equation}
That is, the ($r{>}0$ `sheet' of the) ground state variety is the
conifold
$\shM$, with an exocurve, $A_j$, attached at each singular point.

In the other, $r{<}0$ `sheet' of the field space $\cF$, the first term 
in
the union~(\ref{e:QUnion}) turns out to be empty. On the other hand, 
the
second one is not since $\{p,~s=\sh_j\}$ does include the complex
$p$-plane in which
\Eq{e:gaugeD} shows that $r<0$. In this case, the second term
in the union~(\ref{e:QUnion}) again turns out to be of
the form~(\ref{e:Antennae}), except this time the $A_j$'s have a single
common point, the Landau-Ginzburg orbifold~(\ref{e:LGO}).

Alternatively, consider the symplectic quotient: impose the vanishing of
$D_r$, \ie, intersect with $D_r^{-1}(0)$, and pass to the 
$S^1$-quotient.
To this end, consider each term in the union~(\ref{e:Union}) 
separately.
\topic{$r>0$}
Now $D_r=0$ implies that
\begin{equation}
   \|s\|^2-5|p|^2~=~r~>~0~,\quad\To\quad \|s\|^2~\geq~|r|~>~0~.
    \Label{e:snot0}
\end{equation}
The ground state variety now is the $S^1$-quotient of the union:
\begin{equation}
   \Big\{p{=}0,s\neq\sh_j:G(s)=0,\|s\|^2=r\Big\}
   ~\cup~
    \Big\{p,~s{=}\sh_j:\|s\|^2=r_+\Big\}~, \Label{e:UnionQ}
\end{equation}
where
\begin{equation}
   r_+~=~5|p|^2+|r|~. \Label{e:rplus}
\end{equation}
Note that the $p=0$ points of the second component, where
$\|s\|^2=r_+{=}r$, the
\begin{equation}
   \Big\{(p,s){=}(0,\sh_j):\|s\|^2=r\Big\}\Label{e:Sings}
\end{equation}
points are the $s\to\sh_j$ limiting points of the first component, 
since
$G(\sh_j){=}0$.
  The $S^1$ quotient of these are the (nodal) singular points of the
conifold~(\ref{e:shM}), and they connect the two terms in the
union~(\ref{e:UnionQ}). This then becomes $\shM\cup\sqcup_jA_j$, just 
as
obtained using the holomorphic quotient~(\ref{e:AntennaM}).

\topic{$r<0$}
Now $D_r=0$ implies that
\begin{equation}
   \|s\|^2-5|p|^2~=~r~<~0~,\quad\To\quad |p|~\geq~\sqrt{|r|/5}~>~0~.
    \Label{e:pnot0}
\end{equation}
This renders the first term in the union~(\ref{e:qunion}) empty, and 
the
ground state variety now is:
\begin{equation}
    \Big\{p,~s{=}\sh_j:\|s\|^2=r_-\Big\}\Big/S^1~, \Label{e:LGOQ}
\end{equation}
where now
\begin{equation}
   r_-~=~\|\sh_j\|^2+|r|~. \Label{e:rminus}
\end{equation}
Note that at $s=0=r_-$, where $|p|=\sqrt{|r|/5}$, the
\begin{equation}
   \Big\{p,~s{=}0:5|p|^2=r\Big\}\Big/S^1
\end{equation}
  point is common to all components of the second component, and
is the `fuzzy point' of the Landau-Ginzburg orbifold~(\ref{e:LGO}).

The foregoing proves the following
\begin{lemma} \Label{l:GndStVar}
   With the `ingredients,' $\{(p,s),\hat\lm,D_r,W\}$, defined as in
\Eqs{e:ProjMap}, (\ref{e:gaugeD}) and~(\ref{e:W=pG}), the $r{>}0$
`sheet' of the ground state variety (Definition~\ref{d:GndStVar}),
$\cV^+$, becomes
a stratified variety, $\shM\,\cup\,\sqcup_jA_j$ when $G(s)$ is
non-transversal as specified in \Eq{e:nodalG}.\\
\end{lemma}
\Remark
The `main' stratum~(\ref{e:nonConifold}) has complex dimension 3, while
the `exocurves'~(\ref{e:Antennae}) minus the singular points 
$\shx_j$
form the complex dimension 1 stratum; the singular points,
$\sqcup_j\shx_j=\hbox{Sing}(\shM)$, form the complex dimension 0 
stratum.

\begin{corl}
Under the same conditions as in Lemma~\ref{l:GndStVar}, the $r{<}0$ 
`sheet' of the ground state variety, $\cV^-$, is the stratified variety:
the union of the exocurves~(\ref{e:Antennae}),
$\cup_jA_j$, connected at the `fuzzy point'
of the Landau-Ginzburg orbifold~(\ref{e:LGO}).
\end{corl}
\Remark
The $r{<}0$ stratified variety consists of the exocurves $A_j$ minus
the `fuzzy point' which form the complex dimension 1 stratum, and the 
`fuzzy
point'~(\ref{e:LGO}) which forms the complex dimension 0 stratum.

The resulting non-transversal ground state variety is illustrated in
Fig.~\ref{f:singular}
\begin{figure}[ht]
  \setlength{\unitlength}{.85mm}
  \begin{picture}(170,110)(-20,0)
        \put(0,0){\epsfxsize=140mm\epsfbox{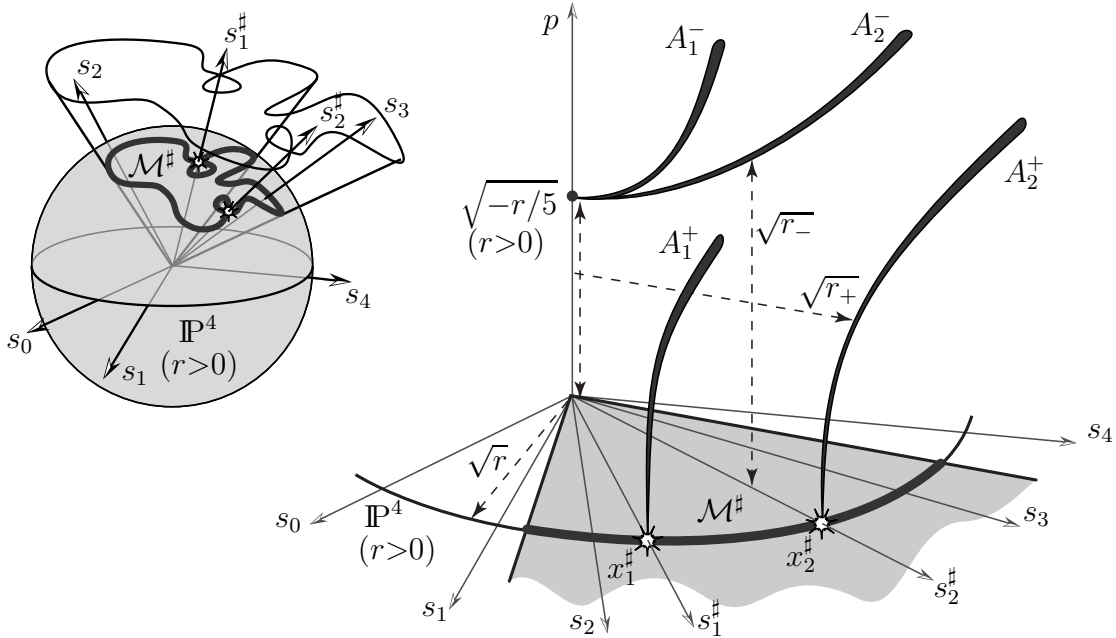}}
        \put(39,17){$s_0$}
        \put(62,3){$s_1$}
        \put(85,1){$s_2$}
        \put(156,18){$s_3$}
        \put(166,31){$s_4$}
        \put(105,1){$\sh_1$}
        \put(142,6){$\sh_2$}
        \put(91,9){$\shx_1$}
        \put(119,11){$\shx_2$}
        \put(53,17){$\IP^4$}
        \put(52,12){$(r{>}0)$}
        \put(69,26){$\sqrt{r}$}
        \put(105,18){$\shM$}
        \put(81,95){$p$}
        \put(68,66){$\sqrt{-r/5}$}
        \put(69,60){$(r{>}0)$}
        \put(114,63){$\sqrt{r_-}$}
        \put(100,92){$A^-_1$}
        \put(129,94){$A^-_2$}
        \put(121,53){$\sqrt{r_+}$}
        \put(99,60){$A^+_1$}
        \put(153,71){$A^+_2$}
        \put(-3,45){$s_0$}
        \put(15,40){$s_1$}
        \put(8,88){$s_2$}
        \put(56,82){$s_3$}
        \put(50,52){$s_4$}
        \put(31,93){$\sh_1$}
        \put(46,80.5){$\sh_2$}
        \put(24,46){$\IP^4$}
        \put(21,41){$(r{>}0)$}
        \put(16,72){$\shM$}
 \end{picture}
  \caption{A non-transversal ground state variety, $\cV\in{\cal F}$,
and its `geometric' phase, $\shM=G^{-1}(0)\cap\IP^4$, in the top
left inset. The rays $\sh_j$ pass through the nodesof $\shM$,
$\shx_j$, which is where the exocurves, $A^+_j$,
attach to $\shM$ in the $r{>}0$ `sheet.'
In the $r{<}0$  `sheet,' the exocurves $A^-_j$ all meet at the
Landau-Ginzburg `fuzzy point.'}
\Label{f:singular}
\end{figure}
\Remark
Since the 'fuzzy point' of the Landau-Ginzburg orbifold~(\ref{e:LGO})
may be, formally, considered as the (negative size) collapse (or,
perhaps more properly, analytic continuation) of the
3-dimensional Calabi-Yau variety $\cM$, the same relation remains
between $\cV^+$ and $\cV^-$, regardless of the (non)transversality
of $G(s)$.

\subsection{The exocurves}
We now turn to study the exocurves, $A_j$, in some detail. In 
particular,
we prove:
\begin{lemma}\Label{l:Antennae+}
In the $r>0$ `sheet' of the field space, $\cF$, the
exocurves~(\ref{e:Antennae}) are
\begin{equation}
   A^+_j~\simeq~\IC\IP^1_{[-5,1]}~\simeq~\IC^1~.
\end{equation}
\end{lemma}
\Proof
In the $r>0$ `sheet,' the definition~(\ref{e:Antennae}) of the 
exocurve:
\begin{equation}
    A^+_j~\define~\{p,\sh_j\}\Big/\hat{\lm}~, \Label{e:Antenna}
\end{equation}
includes implicitly that $\|s\|^2\geq|r|>0$ owing to 
\Eq{e:snot0},
and the superscript `$+$' reminds that $r{>}0$. That is,
\begin{equation}
   (p,\sh_j)~\cong~(\lm^{-5}p,\lm\sh_j)~,\qquad \lm\in\IC^*~, 
\Label{e:proj}
\end{equation}
which defines $A^+_j$ as the weighted projective space
$A^+_j=\IP^1_{[-5,1]}$, proving the first part of~(\ref{e:Antenna}).
  This case, however, differs from the usual consideration
of weighted projective spaces~\cite{dim:weight} in that the weights, 
$-5$
and $1$, are of opposite sign. Still, we proceed by considering the
two candidate charts:
\begin{equation}
   U_p=(p,\sh_j)_p~\cong~(1,u_p)~,\quad p\neq0~,\qquad
    u_p\define\sh_j\,p^{1/5}~,  \Label{e:Up}
\end{equation}
and
\begin{equation}
   U_s=(p,\sh_j)_s~\cong~(u_s,1)~,\quad \sh_j\neq0~,\qquad
    u_s\define p\,(\sh_j)^5~.  \Label{e:Us}
\end{equation}
In both cases, the equivalences are obtained using the 
map~(\ref{e:proj}),
however with $\lm=p^{1/5}$ in the first case, and $\lm=(\sh_j)^{-1}$ in 
the
second. Now, in the second candidate chart, $U_s$, the limit point
$p,u_s\to0$ is included, and so
\begin{equation}
   U_s=(p,\sh_j)_s~\cong~(u_s,1)~\simeq~\IC^1
\end{equation}
is a proper chart. On the other hand, in the first candidate chart, 
$U_p$,
the limit point $\sh_j,u_p\to0$ is {\it excluded\/} by the
inequality~(\ref{e:snot0}), so that
\begin{equation}
   U_p=(p,\sh_j)_p~\cong~(1,u_p)~\simeq~\IC^*
\end{equation}
is {\it not\/} a proper chart. In its place, we should introduce two
$\IC^1$-like charts which cover $U_p$. However, this will not
really be necessary since \Eqs{e:Up} and~(\ref{e:Us}) imply that
\begin{equation}
   u_s=p(\sh_j)^5~\mapsto~u_p^5~, \Label{e:glue}
\end{equation}
which is a 1-to-5 holomorphic map outside $u_s{=}0$. That is,
$U_p=(1,u_p)\simeq\IC^*$ is a 5-fold cover of
$(U_s{-}0)=(u_s,1)_{u_s{\neq}0}\simeq\IC^*$; $u_s{=}0$ is of course the
branching point of this holomorphic covering. Therefore, $A^+_j$ may be
parametrized by $u_s$ and so $A^+_j\simeq U_s$. With~(\ref{e:glue}) as 
the
`glueing map,' we then have that the $j^{th}$ exocurve is:
\begin{equation}
   A^+_j~\simeq~\IP^1_{[-5,1]}~=~U_p \cup U_s ~=~ U_s ~\simeq~ \IC^1.
   \Label{e:lemma2a}
\end{equation}
\QED

\begin{lemma}\Label{l:Antennae-}
In the $r<0$ `sheet' of the field space, $\cF$, the
exocurves~(\ref{e:Antennae}) are
\begin{equation}
   A^-_j~\simeq~\IC\IP^{1,-}_{[-5,1]}~\simeq~\IC^1/\ZZ_5~.
\end{equation}
\end{lemma}
\Proof
In the $r<0$ `sheet,' the definition~(\ref{e:Antenna}--\ref{e:proj}) 
still
guarantees that $A^-_j\simeq\IP^{1,-}_{[-5,1]}$, but now
\Eq{e:pnot0} enforces $|p|^2\geq|r|>0$, as indicated by the
superscript `$-$'.
  We again proceed by considering the two candidate charts~(\ref{e:Up})
and~(\ref{e:Us}). This time, it is in the first candidate chart, $U_p$,
that the limit point $\sh_j,u_p\to0$ is included, and so
\begin{equation}
   U_p=(p,\sh_j)_p~\cong~(1,u_p)~\simeq~\IC^1
\end{equation}
is a proper chart. Similarly, in is now the second candidate chart, 
$U_s$,
from which the limit point $p,u_s\to0$ is {\it excluded\/} by the
inequality~(\ref{e:pnot0}), so that
\begin{equation}
   U_s=(p,\sh_j)_s~\cong~(u_s,1)~\simeq~\IC^*
\end{equation}
is {\it not\/} a proper chart. Again, it is not necessary to introduce 
two
$\IC^1$-like charts to cover $U_s$, since \Eq{e:Up}
and~(\ref{e:Us}) again imply the 1-to-5 holomorphic map~(\ref{e:glue})
now outside $u_p{=}0$. Now $(U_p{-}0)=(1,u_p{\neq}0)\simeq\IC^*$ is a
5-fold cover of $U_s=(u_s,1)\simeq\IC^*$, and $u_p{=}0$ is of course 
the
branching point of this holomorphic covering. Therefore, $A^-_j$ now 
must be
parametrized by $u_p$ which is 5-fold redundant except at $u_p{=}0$.
Therefore, we now have that the $j^{th}$ exocurve is:
\begin{equation}
   A^-_j~\simeq~\IP^{1,-}_{[-5,1]}~=~(U_p \cup U_s)/\ZZ_5~=~ U_p/\ZZ_5
       ~\simeq~ \IC^1/\ZZ_5.  \Label{e:lemma2b}
\end{equation}
\QED

\Remark
Note that the $\ZZ_5$ quotient in Lemma~\ref{l:Antennae-} precisely
corresponds to the $\ZZ_5$ quotient in \Eq{e:LGO}. Indeed, this 
says
that the ``fuzzy point'' of the Landau-Ginzburg orbifold~(\ref{e:LGO})
becomes
\begin{equation}
   (\sqcup_{j=1}^n A^-_j\Big) ~\simeq~ (\IC^1/\ZZ_5)^{\vee5}~,
\end{equation}
\ie, the `plum product' of five copies of the $\IC^1/\ZZ_5$ cone, all
connected at the vertex---the ``fuzzy point''~(\ref{e:LGO}).

The exocurves are illustrated in Fig.~\ref{f:antenna}.
\begin{figure}[ht]
  \setlength{\unitlength}{.85mm}
  \begin{picture}(170,135)(-20,0)
        \put(0,0){\epsfxsize=135mm\epsfbox{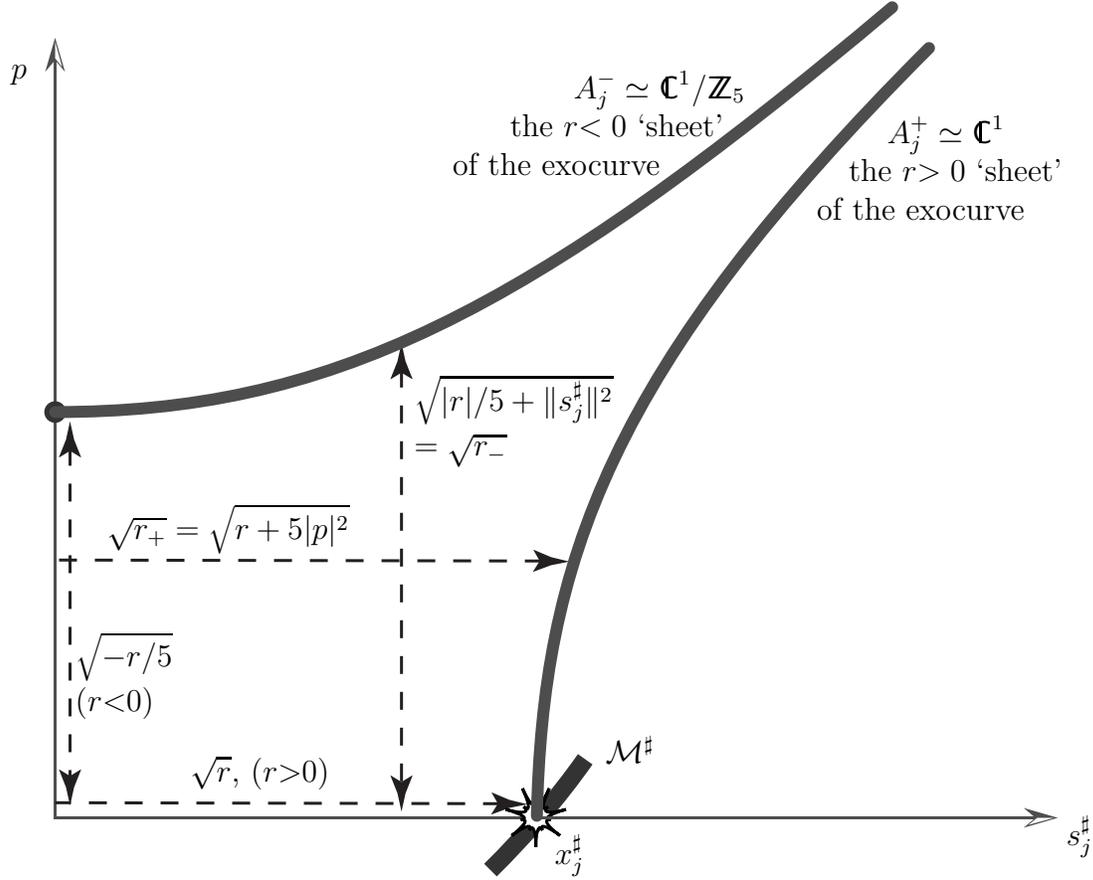}}
        \put(160,5){$\sh_j$}
        \put(80,2){$\shx_j$}
        \put(88,18){$\shM$}
        \put(5,33){$\sqrt{-r/5}$}
        \put(5,26){$(r{<}0)$}
        \put(-5,125){$p$}
        \put(23,15){$\sqrt{r}$, $(r{>}0)$}
        \put(10,53){$\sqrt{r_+}=\sqrt{r+5|p|^2}$}
        \put(132,115){$A^+_j\simeq\IC^1$}
        \put(126,109){the $r{>0}$ `sheet'}
        \put(121,103){of the exocurve}
        \put(58,66){$=\sqrt{r_-}$}
        \put(58,73){$\sqrt{|r|/5+\|\sh_j\|^2}$}
        \put(83,122){$A^-_j\simeq\IC^1/\ZZ_5$}
        \put(73,116){the $r{<0}$ `sheet'}
        \put(64,110){of the exocurve}
  \end{picture}
  \caption{Over the non-transversal rays, $\sh_j$, of $G(s)$ the complex
variable $p$ is subject only to the projectivization action, $\hat\lm$.
The resulting space, the exocurve $\IP^1_{[-5,1]}=\{p,\sh_j\}/\hat\lm$,
is illustrated here for both  the $r{>}0$ `sheet' ($A^+_j$),
and the $r{<}0$ `sheet' ($A^-_j$).}
\Label{f:antenna}
\end{figure}

\subsection{A comparison}
For comparison, we include a similar analysis of $\IP^1_{[5,1]}$. In
contrast to the non-compact $\IP^1_{[-5,1]}$, the weighted projective 
space
$\IP^1_{[5,1]}$ will prove to be compact.

Again, it is possible to view $\IP^1_{[5,1]}$ both as a holomorphic
quotient,
\begin{equation}
   \IP^1_{[5,1]} ~\simeq~ \left\{q,s\right\}/\hat{\m}
\end{equation}
where
\begin{equation}
   \hat\m\colon (q,s) \mapsto (\m^5q,\m s)~,\qquad \m\in\IC^*~,
    \Label{e:projmap}
\end{equation}
and also as a symplectic quotient,
\begin{equation}
   \IP^1_{[5,1]} ~\simeq~
    \left\{\left\{q,s\right\}\cap \D_r^{-1}(0)\right\}/S^1
\end{equation}
where the $S^1$-action is the restriction of~(\ref{e:projmap}) to 
$|\m|=1$
and \Eq{e:gaugeD} now becomes
\begin{equation}
   \D_r~=~|s|^2+5|q|^2-r~,\qquad r\in\IR~. \Label{e:gauged}
\end{equation}
Following the proofs of Lemmas~\ref{l:Antennae+} and~\ref{l:Antennae-}, 
we
consider the latter.

The vanishing of $\D_r$ now simply states 
that
\begin{equation}
   |s|^2+5|q|^2 ~=~r~\geq~0~, \Label{e:rPos}
\end{equation}
where $r=0$ would force $s=0=q$, the trivial solution. Restricting
then to $r\neq0$, \Eq{e:rPos}
implies the positive definiteness of $r$ on $\D_r^{-1}(0)$,
so that there is only the $r>0$ sheet. Indeed, this is
precisely why the $r<0$ sheet of $\IP^1_{[-5,1]}$ appears to be rather
unfamiliar an object.
Owing to the inequality~(\ref{e:rPos}),
$s,q$ must not vanish {\it simultaneously\/}; either one of them
however may very well vanish while the other one is nonzero. Thus, 
unlike
in the case of $\IP^1_{[-5,1]}$, we now have two perfectly proper 
coordinate
charts:
\begin{equation}
   U_q=(q,s)_q~\cong~(1,u_q)~\simeq~\IC^1~,\quad\hbox{using}\quad
    \m=q^{-1/5}~,~ q\neq0~,\quad u_q=s\,q^{-1/5}~,
\end{equation}
and
\begin{equation}
   U_s=(q,s)_s~\cong~(u_s,1)~\simeq~\IC^1\quad\hbox{using}\quad
    \m=s^{-1}~,~s\neq0~,\quad u_s=q\,s^{-5}~.
\end{equation}
The two chart coordinates, $u_q$ and $u_s$, respectively, {\it can\/}
attain the value of 0, since $s,u_q\to0$ is permitted in $U_q$ where
$q\neq0$, and $q,u_s\to0$ is permitted in $U_s$ where $s\neq0$. 
Finally,
the two charts are glued through the relation
\begin{equation}
   u_s~=~u_q^{-5}~,\qquad\hbox{where}\quad u_s,u_q\neq0~, 
\Label{e:glue1}
\end{equation}
which provides a 1-5 map:
\begin{equation}
   \left\{U_s-0\right\} ~\tooo{~1-5~}~ \left\{U_q-0\right\}~.
\end{equation}
To render the map~(\ref{e:glue1}) single-valued, we may glue together 
$U_s$
and $U_q/\ZZ_5$: $0\in U_s$ becomes `$\infty$' added to $U_q/\ZZ_5$, and 
$0\in
U_q/\ZZ_5$ becomes `$\infty$' added to $U_s$. The resulting space,
$\IP^1_{[5,1]}=U_s\cup(U_q/\ZZ_5)$, then is compact and smooth except 
at
$0\in U_q/\ZZ_5$, where $\IP^1_{[5,1]}$ has a $\ZZ_5$ quotient 
singularity,
\ie, a conical singularity with $(1-\inv5)2\p=8\p/5$ deficit angle.

Note, however, that by Delorme's Lemma~\cite{dol:weight}, $\IP^1_{[k,1]}\approx\IP^1_{[1,1]}\equiv\IP^1\simeq S^2$.
The relationship `$\approx$' here denotes a $k$-to-$1$ map of the
coordinates as used here, but an isomorphism of the corresponding
coordinate rings, which then extends to an isomorphism of the
respective spaces~\cite{dol:privco}.

\subsection{A one-point compactification of exocurves}
For later convenience and use, we describe here a one-point
compactifications of the exocurves.

It is straightforward from the proof of
Lemmas~(\ref{l:Antennae+}) that the limiting point
$s,u_p\to\infty$ may be added to $U_p$. Upon the inversion of its
variables, this now becomes a proper coordinate chart:
\begin{equation}
   \Tw{U}_p = (1,w_p)~\simeq~\IC^1~,\qquad
    w_p\define u_p^{-1}={\sh}_j^{-1}p^{-1/5}~.
\end{equation}
Clearly, the glueing map now becomes
\begin{equation}
   u_s~=~w_p^{-5}~\colon~
   \left\{U_s-0\right\}~\tooo{~1-5~}~\left\{\Tw{U}_p-0\right\}~.
    \Label{e:glue2}
\end{equation}
To render the glueing map~(\ref{e:glue2}) single-valued, we form
\begin{equation}
 \bar{A}{}_j^+~\define~U_s\cup(\Tw{U}_p/\ZZ_5)~\simeq~\IP^1_{[5,1]}~.
  \Label{e:CompAnt}
\end{equation}
As shown above, this is compact and isomorphic to $\IP^1\simeq S^2$.

\section{Cohomology and Homology of $\bar\cV$}
\Label{s:Hom}
As presented in Theorem 1, the stratified variety can be written as
\begin{equation}
   \bar \cV~=~\shM ~\cup~ \bigcup_{j=1}^n \bar {A}_j~, \qquad
   \shM \cap \bar{A}_j = \shx_j~.  \Label{e:AntennaM2}
\end{equation}

The Mayer-Vietoris principle then induces the long exact cohomology sequence
\begin{equation}
 \3\to H^q(\bar\cV) \to H^q(\shM)\oplus H^q(\cup_j\bar{A}_j)\to
         H^q(\shM\cap\cup_j\bar{A}_j)\to H^{q+1}(\bar\cV)\to\3
 \Label{e:LongHomSeq}
\end{equation}
Since
\begin{equation}
 \shM\cap\cup_j\bar{A}_j=\cup_j(\shM\cap\bar{A}_j)=\sqcup_j\shx_j~,
\end{equation}
owing to our assumption that $\shx_j$ are isolated (non overlapping) nodes, and
\begin{equation}
 H^q(\cup_j\shx_j)=\bigoplus_{j=1}^n H^q(\shx_j)\simeq\d_{q,0}\IC^{\oplus n}~,
\end{equation}
the long exact sequence~(\ref{e:LongHomSeq}) breaks into five isomorphisms:
\begin{equation}
 H^q(\bar\cV) = H^q(\shM)\oplus H^q(\cup_j\bar{A}_j)~,
 \qquad\hbox{for~}q=2,{\cdots},6~, \Label{e:IsoMV}
\end{equation}
and
\begin{equation}
 0\to H^0(\bar\cV)\tooo{~\a~} H^0(\shM)\oplus H^0(\cup_j\bar{A}_j)\tooo{~\b~}
 H^0(\sqcup_j\shx_j)\to H^1(\bar\cV)\to0~.
\end{equation}
The above map $\b$ is induced from the injective inclusion $\sqcup_j\shx_j=\shM\cap(\sup_j\bar{A}_j)\to\shM\sqcup(\cup_j\bar{A}_j)$, and so is surjective. Then:
\begin{equation}
 H^1(\bar\cV)=\emptyset~, \quad\hbox{and}\quad
 H^0(\bar\cV) \simeq \IC~. \Label{e:BegMV}
\end{equation}

\subsection{Contributions from the antennae}
Recall from a previous section that $\bar {A_j} \simeq\IP^1\simeq S^2 $. Then,
\begin{equation}
 H^q(\bar{A}_j)~~\left\{
 \begin{array}{ll}
   \simeq\IC^1 &\mbox{for~} q=0,2~,\\[1mm]
   =\emptyset &\mbox{otherwise.}
 \end{array}\right.
\end{equation}
With this, for $\shM$ with $n$ simple nodes, \Eqs{e:IsoMV} and~(\ref{e:BegMV}) would seem to imply that $H^q(\cV)$ should equal to $H^q(\shM)$, except for $q=2$, where it ought to be augmented by
$H^2(\cup_j\bar{A}_j)\simeq\oplus\IC^{\oplus n}$.

This, however, is not correct: the (area) 2-forms of the $n$ antennae are not independent cohomology elements. As described in detail in Ref.~\cite{hub:cymani}, $N$ mutually exclusive subsets of the $n$ $\shx_j$'s lie on corresponding 4-cycles $C^{(4)}_k\subset\shM$, $k=1,{\cdots},N$. Let $J_k$ denote the multiindex containing the indices, $j$, of all $\shx_j$'s that lie on $C^{(4)}_k$.
 Clearly then,
\begin{equation}
 \bar{A}_j \cap C^{(4)}_k ~=~ \left\{
  \begin{array}{ll}
   \shx_j & \mbox{if~}j\in J_k\mbox{~{\it i.e.}~}\shx_j\in C^{(4)}_k~,\\[1mm]
   \emptyset & \mbox{otherwise.}
  \end{array}\right. \Label{e:ClasSing}
\end{equation}
Considering then the {\it homology\/} elements in $H_q(\bar\cV)$, dual to the cohomology group obtained in \Eqs{e:IsoMV} and~(\ref{e:BegMV}), and denoting them by square brackets, we have:
\begin{equation}
 [\bar{A}_j] \cap [C^{(4)}_k] ~=~ \left\{
  \begin{array}{ll}
   1 & \mbox{if~}j\in J_k~,\\[1mm]
   0 & \mbox{otherwise.}
  \end{array}\right. \Label{e:DualHom}
\end{equation}
Owing to this result, it follows that
\begin{equation}
 \left.\begin{array}{lll}
  [\bar{A}_j]=[\bar{A}_{j'}] &\mbox{if~~}j,j'\in J_k~,\\[1mm]
  [\bar{A}_j]\neq[\bar{A}_{j'}] &\mbox{otherwise.}
       \end{array}\right\} \Label{e:EquAnt}
\end{equation}
That is, the $n$ antennae, $\{\,\bar{A}_j\,\}$ contribute only $N$ inequivalent 2-cycles, so
\begin{equation}
 H_2(\cup_j\bar{A}_j) ~\simeq~ \IC^{\oplus N} ~\simeq~ H^2(\cup_j\bar{A}_j)~.
 \Label{e:AntCoh}
\end{equation}

\subsection{The combined result}
Combining \Eqs{e:IsoMV}, (\ref{e:BegMV}) and~(\ref{e:AntCoh}) proves
\begin{lemma}\Label{l:Coh}
Let $\bar\cV$ as defined in \Eq{e:AntennaM2}, where $\shM$ is a conifold with only $n$ isolated nodes ($\shx_j$) lying on $N$ distinct 4-cycles $C^{(4)}_k$, and $\bar{A}_j$ as defined in \Eq{e:CompAnt}. Then,
\begin{equation}
 H^q(\bar\cV) = \left\{
  \begin{array}{l}
   H^q(\shM)\quad \mbox{for~} q\neq2~,\\[1mm]
   H^2(\shM)\oplus H^2(\cup_j\bar{A}_j) \simeq H^2(\shM)\oplus\IC^{\oplus N}~.
  \end{array}\right.
\end{equation}
\end{lemma}

As it stands, with $H^3(\bar\cV)=H^3(\shM)$, the {\it complete} $H^*(\bar\cV)$ can have neither Poincar\'e duality nor a Hodge decomposition. Both are obstructed by the fact that the 3-cycle(s) which pass through the $\shx_j$'s remain without dual 3-cycle(s)~\cite{hub:cymani}. In fact, the subgroup of $H^3(\bar{\cV})$ generated by the 3-cycles passing through the $\shx_j$'s may well be odd-dimensional, making this obstruction manifest.

However, $\oplus_qH^{2q}(\bar\cV)$ subgroup does exhibit both Poincar\'e duality and a Hodge decomposition. As usual, $H^2(\bar{A}_j)\simeq\IC^1$ is generated by the volume $(1,1)$-form on $\bar{A}_j\simeq\IP^1$. Moreover, dually to the homology result~(\ref{e:DualHom}), the  volume $(1,1)$-forms, $\w_{(1,1)}^j$, of all $\bar{A}_j$'s which intersect $C^{(4)}_k$ are dual to the $(2,2)$-form $\w_{(2,2)}^k$, itself dual to $C^{(4)}_k$. In fact, the double dualities\footnote{By $\dual$ we denote the standard homology--cohomology duality, and use $\Dual$ for Poincar\'e duality.}
\begin{equation}
 \w_{(1,1)}^j\dual[\bar{A}_j]\Dual[C^{(4)}_k]~, \quad\mbox{and}\quad
 [\bar{A}_j]\Dual[C^{(4)}_k]\dual\w_{(2,2)}^k~, \qquad\mbox{for~}j\in J_k~,
 \Label{e:DDual}
\end{equation}
establishes the isomorphisms
\begin{equation}
 \w_{(1,1)}^j\simeq[C^{(4)}_k]~, \quad\mbox{and}\quad
 [\bar{A}_j]\simeq\w_{(2,2)}^k~, \qquad\mbox{for~}j\in J_k~,
 \Label{e:Isom}
\end{equation}
whereupon \Eq{e:EquAnt} implies that also
\begin{equation}
 \left.\begin{array}{lll}
  [\w_{(1,1)}^j]=[\w_{(1,1)}^{j'}] &\mbox{if~~}j,j'\in J_k~,\\[1mm]
  [\w_{(1,1)}^j]\neq[\w_{(1,1)}^{j'}] &\mbox{otherwise.}
       \end{array}\right\} \Label{e:EqVAnt}
\end{equation}
On a more fomal level, \Eq{e:DualHom} implies that also
\begin{equation}
 [\w_{(1,1)}^j] \cup [\w_{(2,2)}^k] ~=~ \left\{
  \begin{array}{ll}
   1 & \mbox{if~}j\in J_k~,\\[1mm]
   0 & \mbox{otherwise.}
  \end{array}\right. \Label{e:DualCoh}
\end{equation}
Clearly, the evaluation map of the cup product here cannot be the integration (of the wedge product of the indicated forms) over the stratified variety $\bar\cV$ in any conventional sense. Instead, it may be taken to reduce to the evaluation over the point of common support, $\shx_j=\bar{A}_j\cap C^{(4)}_k$ if $j\in J_k$, and is vacuous otherwise.

Owing to the isomorphisms~(\ref{e:Isom}), the $[\bar{A}_j]\Dual[C^{(4)}_k]$ duality (when $j\in J_k$) implies the desired Poincar\'e duality of $\w_{(1,1)}^j\Dual\w_{(2,2)}^k$, for all $j\in J_k$. Let $\vev{\w_{(2,2)}^k}$ denote the subgroup of $H^{(2,2)}(\shM)$ generated by the $\w_{(2,2)}^k$'s. The quotient $H^{(2,2)}(\shM)/\w_{(2,2)}^k$ is then generated by the $(2,2)$-forms dual to 4-cycles which do not pass through $\shx_j$; this quotient is easily seen to form an additive group, exhibiting both poincar\'e duality and Hodge decomposition.

The foregoing then proves:
\begin{lemma}\Label{l:CohPH}
Let $\bar\cV$ as defined in \Eq{e:AntennaM2}, where $\shM$ is a conifold with only $n$ isolated nodes ($\shx_j$) lying on $N$ distinct 4-cycles $C^{(4)}_k$, and $\bar{A}_j$ as defined in \Eq{e:CompAnt}. Then,
\begin{equation}
 \oplus_q H^{2q}(\bar\cV) = \left\{
  \begin{array}{l}
   H^{2q}(\shM)\quad \mbox{for~} q\neq1~,\\[1mm]
   H^2(\shM)\oplus H^2(\cup_j\bar{A}_j) \simeq H^2(\shM)\oplus\IC^{\oplus N}~.
  \end{array}\right.
\end{equation}
has both an induced Hodge decomposition and Poincar\'e duality, as induced by the double dualities~(\ref{e:DDual}).
\end{lemma}

\section{Deformations, Resolutions and the Mirror Map}
\Label{s:DefRes}
We have originally restricted $\shM$ to conifolds with only nodes (\ie, double points, or $A_1$ hypersurface singularities), $\shx_j$.
Their local neighborhood is isomorphic to the cone $\IC^4/Q$, where $Q$ is a non-degenerate quadratic polynomial over $\IC^4$.
In a small resolution, this neighborhood is replaced with a copy of the total space of an ${\cal O}(-1,-1)\define{\cal O}(-1){\oplus}{\cal O}(-1)$
bundle over $\IP^1\simeq S^2$.
 In short, a small resolution replaces each node $\shx_j$ with a $(-1,-1)$-curve, $\IP^1_{\natural,j}\simeq S^2$. 
 Since there are two topologically distinct ways to do this at each node $\shx_j$, a conifold $\shM$ with $n$ would appear to
have $2^n$ small resolutions, $\nM_I$. However, all nodes $\shx_j$ which lie on a single 4-cycle $C^{(4)}_k\subset\shM$ must be
resolved ``compatibly'': all the corresponding 2-spheres $\IP^1_{\natural,j}\in\nM_I$ intersect $C^{(4)}_{\natural,k}\in\nM_I$
(the proper transform of $C^{(4)}_k\subset\shM$) in a single point and so must all represent the same element of $H_2(\nM_I)$,
the one that is dual to $C^{(4)}_{\natural,k}$. With the use of \Eq{e:ClasSing}, this implies that
\begin{equation}
 \left.\begin{array}{ll}
  [\w_{(1,1)}^{\natural,j}]=[\w_{(1,1)}^{\natural,j'}] &\mbox{if~~}j,j'\in J_k~,\\[1mm]
  [\w_{(1,1)}^{\natural,j}]\neq[\w_{(1,1)}^{\natural,j'}] &\mbox{otherwise,}
       \end{array}\right\} \Label{e:EqVAntSmR}
\end{equation}
and
\begin{equation}
 [\w_{(1,1)}^{\natural,j}] \cup [\w_{(2,2)}^{\natural,k}] ~=~ \left\{
  \begin{array}{ll}
   1 & \mbox{if~}j\in J_k~,\\[1mm]
   0 & \mbox{otherwise.}
  \end{array}\right. \Label{e:DualCohSmR}
\end{equation}
for
\begin{eqnarray}
 H_2(\nM_I)\ni\IP^1_{\natural,j}
 &\dual& \w_{(1,1)}^{\natural,j}\in H^{(1,1)}(\nM_I)~, \\
 H_4(\nM_I)\ni C^{(4)}_{\natural,k}
 &\dual& \w_{(2,2)}^{\natural,k}\in H^{(2,2)}(\nM_I)~.
\end{eqnarray}
Of course, in \Eq{e:DualCohSmR}, the cup product is indeed obtained as the ordinary wedge product, integrated over the (smooth) manifold $\nM_I$
Consequently, the multiplicity of small resolutions to $I=1,{\cdots},2^N$, where $N$ is the number of $H_2(\nM_I)$ elements which the small
resolution exceptional sets, $\IP^1_{\natural,j}$ represent, \ie, the number of $H_4(\nM_I)$ elements, $C^{(4)}_{\natural,k}$, which are the
proper transforms of the 4-cycles that pass through the nodes $\shx_j\in\shM$.~\cite{hub:cymani}.

The formal identity of the \Eqs{e:EqVAnt}--(\ref{e:DualCoh}) with the \Eqs{e:EqVAntSmR}--(\ref{e:DualCohSmR}) then proves:
\begin{lemma} \Label{l:TriFlop}
   Let $\shM$ be a Calabi-Yau complex 3-dimensional algebraic variety with only a finite number of isolated nodes, $\shx_j$. Let $\nM_1$ 
and $\nM_2$ denote two small resolutions of $\shM$, related by a flop: $\nM_1\buildrel{{\rm f}}\over{\longleftrightarrow}\nM_2$.
Finally, let $\bar{\cV}$ be the compactification of the stratified variety~(\ref{e:AntennaM2}). Then the flop
involution $\nM_1\buildrel{{\rm f}}\over{\longleftrightarrow}\nM_2$ generalizes to a triple of (exo)flops:
\begin{equation}
  \mbox{\begin{picture}(100,45)(0,0)
        \put(-5,0){\epsfxsize=100mm\epsfbox{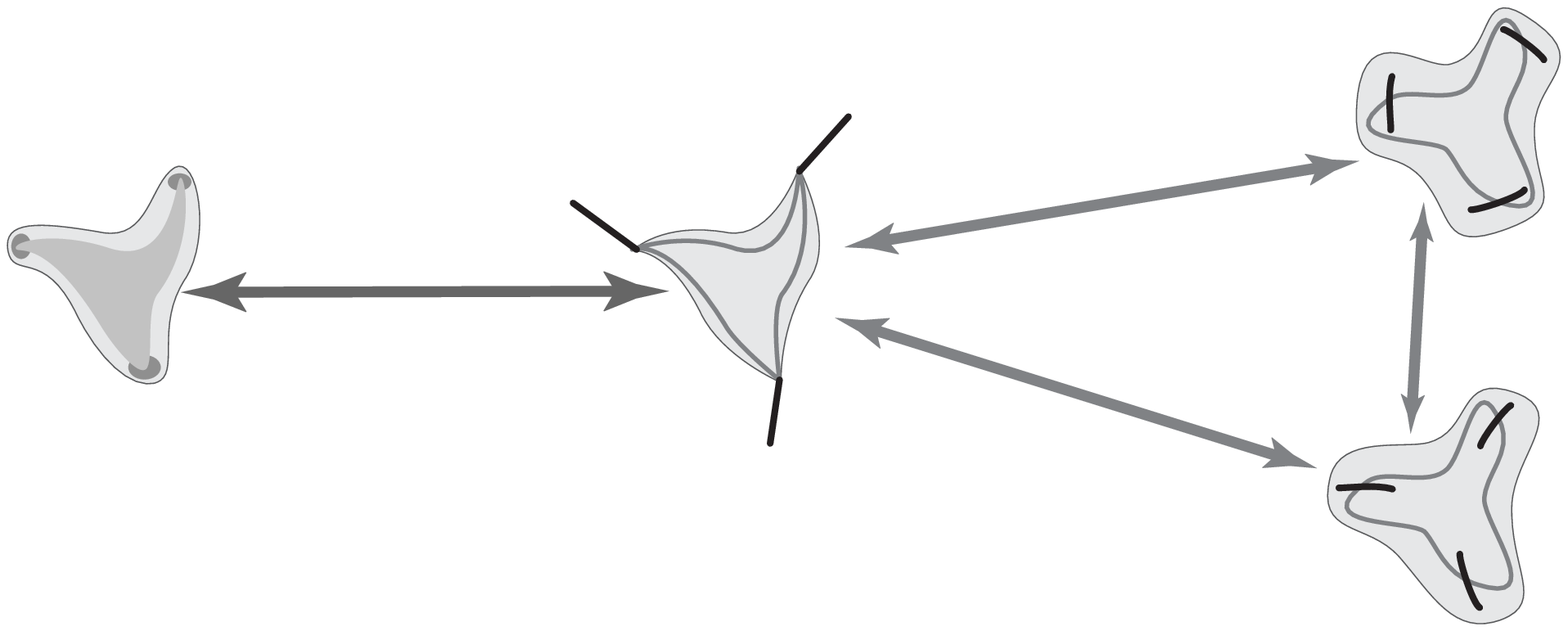}}
        \put(-12,25){$\cM^\flat$}
        \put(16,24){defo}
        \put(38,30){$\bar{\cV}$}
        \put(58,31){exoflop}
        \put(52,10){exoflop}
        \put(95,33){$\nM_1$}
        \put(95,7){$\nM_2$}
        \put(91,19){flop}
  \end{picture}}
\end{equation}
\end{lemma}
The map ``{\it defo}'' to the left is realized as follows: Deformations smooth $\shM$ by replacing the local cones $\IC^4/Q$ centered at each node, $\shx_j$, with a real 3-bundle over a copy of $S^3$. It is easy to see that a deformation of $G(s)$ from the non-transversal choice~\Eq{e:nodalG} to a transversal choice in section~\ref{s:TransV} precisely induces the smoothing of the ground state variety from a (compactified) stratified variety of the type described in section~\ref{s:SingV} to a smooth Calabi-Yau 3-fold of the type described in section~\ref{s:TransV}. This provides the map
$\cM^\flat\buildrel{{\rm defo}}\over{\longleftrightarrow}\bar{\cV}$ in the diagram in Lemma~\ref{l:TriFlop}.

Finally, we note that the above described homology of $\bar{\cV}$ excluding the middle dimension, which we have not discussed herein,
satisfies the requirements given in Ref.~\cite{hub:strcoh}, and so is compatible with the `mirror map.' The extension of this result
to include the (co)homology groups in the middle dimension remains an open question for now and we hope to return to it in a future effort.

\paragraph{Acknowledgments:}
\Label{s:Ack}
The authors would like to thank M. Goresky (Institute for Advanced Study, Princeton, NJ) and N. Ramachandran (University of Maryland, 
College Park, MD) for their helpful comments.


\begin{thebibliography}{text} \rightskip=0pt plus1fill

\bibitem{agm:chngtp}
  P.~Aspinwall, B.~Greene and D.~Morrison:
  \NP{B416}(1994)414--480 .

\bibitem{chsw:strcys}
  P.~Candelas, G.~Horowitz, A.~Strominger and E.~Witten: 
  \NP{B258}(1985)46.

\bibitem{cgh:confld}
  P.~Candelas, P.S.~Green and T.~H\"{u}bsch: \NP{B330}(1990)49--102.

\bibitem{dim:weight}
  A.~Dimca: {\it Singularities and Coverings of Weighted Complete
  Intersections\/}, {\it R.\,Ang.\,J.\,Math.\/}
  {\bf366}(1986)184--193.

\bibitem{dol:weight}
  I.~Dolgachev: {\it Weighted Projective Varieties\/} in Group
  Actions and Vector Fields, {\it Lect.\,N.\,Math.}\,{\bf 956}(1982) 
34--71.

\bibitem{dol:privco}
  I.~Dolgachev: private communication.

\bibitem{ful:torvar}
  W.~Fulton: {\it Toric Varieties\/} (Princeton University Press,
  Princeton, 1990).

\bibitem{gmp:strata}
  J. Cheeger, M. Goresky and R. MacPherson, {\itshape Ann. Math. Studies} {\bf 102} (1982) 303-340.

\bibitem{gsw:string}
  M.B.~Green, J.H.~Schwarz and E.~Witten: {\itshape String Theory},
  (Cambridge University Press, Cambridge, 1987).

\bibitem{gun:torvar}
  E.~Gunter: {\it Convex Bodies and Toric Varieties\/} 
  (Springer-Verlag, New York, 1990)

\bibitem{hub:cymani}
  T.~Hubsch: {\itshape Calabi-Yau Manifolds: A Bestiary for 
Physicists},
  (World-Scientific, Singapore, 1992).

\bibitem{hub:strcoh}
  T.~Hubsch; \MPL{A12}(1997)521.

\bibitem{pol:string}
  J.~Polchinski: {\itshape String Theory},
  (Cambridge University Press, Cambridge, 1997).

\bibitem{str:nullBH}
  A.~Strominger: \NP{B451}(1995)96-108.

\bibitem{wit:phases}
  E.~Witten: \NP{B403}(1993)159--222.

\bibitem{yau:mirror}
  S.-T.~Yau, ed.: {\it Mirror Manifolds\/}
  (International Press, Hong-Kong,
  1990);\\
  B.~Greene and S.-T.~Yau, eds.: {\it Mirror Manifolds II\/}
  (International Press, Hong-Kong, 1996).

\end{thebibliography}
\end{document}